\documentclass[reqno]{amsart}
\usepackage{amsthm}
\usepackage{amsmath}
\usepackage{enumerate}
\usepackage{graphicx}
\usepackage{alphabeta}

\newtheorem{thm}{Theorem}
\newtheorem{cor}{Corollary}

\begin{document}

\title{New bounds for numbers of primes in element orders of finite groups}

\author[C. Bellotti]{Chiara  Bellotti}
\address{School of Science\\
UNSW Canberra at ADFA\\
ACT, Australia}
\email{chbellotti@gmail.com}

\author[T.~M. Keller]{Thomas Michael Keller}
\address{Department of Mathematics\\
Texas State University\\
601 University Drive, San Marcos, TX 78666, USA}
\email{keller@txstate.edu}

\author[T.~S. Trudgian]{Timothy S. Trudgian}
\address{School of Science\\
UNSW Canberra at ADFA\\
ACT, Australia}
\email{t.trudgian@adfa.edu.au}

\date{\today}

\begin{abstract}
Let $\rho(n)$ denote the maximal number of different primes that may occur in the order of a finite solvable group $G$, all elements of which have orders divisible by at most $n$ distinct primes. We show that $\rho(n)\leq 5n$ for all $n\geq 1$. As an application, we improve on a recent bound by Hung and Yang for arbitrary finite groups.
\end{abstract}

\maketitle

\section{Introduction and statement of results}
Let $G$ be a finite solvable group, and let $\omega(n)$ denote the number of distinct prime divisors of $n$. We define $\omega(G)$ as the maximum of $\omega(\textrm{ord}(g))$, over all $g\in G$. Then, finally, we define $\rho(n)$ as the maximum of $\omega(|G|)$, taken over all finite solvable groups $G$ with $\omega(G)= n$. Work on bounding $\rho(n)$ has been undertaken by the second author in the series of papers \cite{Keller2,KellerAB,Keller4}: see also \cite{Moreto21} for some related problems. He \cite{Keller} showed there is a function $C(n)$ such that 
\begin{equation}\label{gufo}
\rho(n) \leq C(n) n,
\end{equation} 
where $C(n)$ is bounded, and $\lim_{n\rightarrow\infty} C(n) = 4$. Of interest is determining the smallest constant $K$ such that
\begin{equation*}\label{gnocchi}
\rho(n) \leq Kn, \quad (n\geq 1).
\end{equation*}
Yang \cite{Yang} showed that $K\leq 17/3$. A weaker (but simpler) proof that $K\leq 7$ was given by Hung and Yang \cite{Hung}. It is conjectured that $K=3$.

If one gets a good value of $K$, then, as outlined in Moret\'{o} \cite{Moreto21}, one can deploy this directly to some related problems. For example, Yang and Qian \cite{YQ} examined the analogous problem for character codegrees and obtained a bound of the shape $K+3$. This was improved recently by Yang \cite{Yang} to $K+2$. At almost the same time Moret\'{o} \cite[Corollary 1.4]{Moreto21} reduced this to $K+ 3/2$. Moret\'{o} notes (see the commentary before Theorem 1.3 in \cite{Moreto21}) that this result is sharp. This shows that a lot of interest has been expended on bounds involving $K$. We use a mix of theory and computation to prove the following.
%This bound of $K+3/2$ is given (improved?) to 7 --- see the proof of Theorem 3.3. in \cite{Hung}).
\begin{thm}\label{gelato}
Let $G$ be a finite solvable group. For all $n\geq 1$ we have $\rho(n) \leq 5n$.
\end{thm}
The same methods leading to the proof of Theorem \ref{gelato} allow us to show, for example, that $\rho(n) \leq 4.1n$ for all $n\geq 1000$. However, as will be made clear in the proof of Theorem \ref{gelato}, it is difficult to remove this restriction on $n$. 
%In short: unlike many explicit problems in number theory, it is not possible to defer the evaluation of the bound for $n\leq 1000$ to computation.

As a result of Theorem \ref{gelato} we can obtain the following improvement on the size of $\rho_{g}(n)$ for arbitrary finite groups. We define $\rho_{g}(n)$ in the same way as $\rho(n)$ except that for the former we consider arbitrary finite groups $G$, not necessarily solvable. 
\begin{cor}\label{bistecca}
Let $G$ be a finite group. For all $n\geq 1$ we have $\rho_{g}(n) \leq 141n^{4}$. 
\end{cor}
This improves on work first examined by Zhang \cite{Zhang}, and improved by the second author (see Remark 16.19(b) in \cite[p.\ 213]{Huppert})  and then Moret\'{o} \cite{Moreto05} to a bound of the form $\rho_{g}(n) \leq A n^{4} \log n$ for some unspecified $A$. Hung and Yang \cite{Hung} proved a version of Corollary \ref{bistecca} with 210 in place of 141.

The outline of this paper is as follows. In \S \ref{biscotto} we collect some relevant background material. This allows us to prove Theorem \ref{gelato} in \S \ref{delfino}. We give the proof of Corollary \ref{bistecca} in \S \ref{foca}, and outline some potential future work stemming from Theorem \ref{gelato}.

\section{Background}\label{biscotto}
We begin by recalling the observation by the second author \cite[p.\ 651]{Keller} that
\begin{equation}\label{oboe}
\rho(n) \leq \pi(\rho(n)) + 4n,
\end{equation}
where $\pi(x)$ counts the number of primes not exceeding $x$. Since $\pi(x) \sim x/\log x$, then, given any $\epsilon$ we can find an $n_{0}$ for which $n\geq n_{0}$ guarantees that $\pi(n) \leq (1 + \epsilon)n/\log n$. Hence rearranging (\ref{oboe}) gives
\begin{equation}\label{bassoon}
\rho(n) \leq 4n \left( 1 - \frac{1+ \epsilon}{\log n_{0}}\right)^{-1}.
\end{equation}
Others have, hitherto, used the simple (and elegant) bound due to Rosser and Schoenfeld \cite[Cor.\ 1]{RS}, namely that 
\begin{equation}\label{anatra}
\pi(n) \leq 1.25506 n/\log n, \quad (n\geq 2).
\end{equation}
When $n$ is small, say, $n<10^{8}$ we could alternatively compute $\pi(n)$ exactly on a computer. 

For larger values of $n$, we could use, for example, the result by Buethe \cite{Buethe2} that $\pi(x) \leq \textrm{li}(x)$ for all $2\leq x \leq 10^{19}$, where $\textrm{li}(x) = \int_{2}^{x} dt/\log t$. We can then bound $\textrm{li}(x)$ by integrating-by-parts and truncating the error term, as in Bennett et al.\ \cite{Bennett} or in Saouter et al.\ \cite{Saouter}. This can be augmented with by another result by Buethe \cite{Buethe}, which has been extended recently by Johnston \cite{DJ}, using the partial verification of the Riemann hypothesis by Platt and the third author \cite{PT}. Namely, we have $\pi(x) - \textrm{li}(x) \leq 1/(8\pi) \sqrt{x} \log^{2}x$ for $x\leq 10^{25}$.  Finally, for still larger values of $n$ we could use any of the (many) results in Broadbent et al. \cite{Broadbent}.

The difficulty with the above is that we run into troubles for small $n$. %Some of the material above can still be used to get good bounds, perhaps something like $\rho(n) \leq 4.01n$ for all $n\geq 10^{10}$, or something similar. 
If we want a bound that holds for all $n$ we are limited by the poor bounds on $\rho(n)$ for small $n$, mainly for $n\leq 100$.

Indeed, Yang's bound of $K\leq 17/3$ comes about from the `worst' bound on $\rho(n)$ occurring at $n=9$. We make an improvement below to show that worst cases are now due to $n=8$ and $n=9$, which give us $K\leq 5$.

Having stressed above the need to examine values of $\rho(n)$ for small $n$, we proceed to list what is known.
We have, from \cite{Keller2,Keller4} that 
\begin{equation}\label{clarinet}
\rho(1) = 2, \quad \rho(2) = 5, \quad \rho(3) = 8, \quad \rho(4)= 12.
\end{equation}
We also have the following collection of bounds (see \cite[\S 1.4]{Keller2} and \cite{Zhang} and \cite{Keller4}) 
\begin{equation}\label{flute}
\rho(n+1) \leq \rho(n) + n + 2, \quad \rho(n) \leq \frac{n(n+3)}{2} -2, \quad \rho(n)\geq 2n,
\end{equation}
where the second bound is true for $n\geq 4$, and the others for all $n\geq 1$.
%In case we need it later, here is another fact to recall: Keller proved that $\rho(n) \geq 2n$ for every $n$.

The exact values of $\rho(n)$ in (\ref{clarinet}) and the second formula in (\ref{flute}) give the following pairs, in which $C(n)$ is the bound in (\ref{gufo}):
\begin{equation}\label{method0}
(n, C(n)) = \{ (5, 18), (6, 25), (7, 33), (8, 42), (9, 52), (10, 63)\}.
\end{equation}
Yang makes a small improvement when $n=9$ by using (\ref{oboe}). 
If $\rho(9) = 52$, then, since $\pi(52) = 15$ we have that the right-hand side of (\ref{oboe}) is 51, whence $\rho(9) \leq 51$, a contradiction. Hence $\rho(9) \leq 51$; this method is unable to reduce $\rho(9)$ further, though small improvements are possible for $n\geq 10$. This (and some additional computation to examine some larger values of $n$) leads to Yang's version of our Theorem \ref{gelato} with constant $51/9 = 17/3$ in place of our $5$.

\section{A new recursive argument}\label{delfino}%\footnote{\textbf{From TST:} May need some help from Thomas here for the theory. And by `some', I mean `much'!}
For a given $n\geq 5$, let $\rho(n) = k$. Let $G$ be a finite solvable group such that $\omega(G)=n$ and $\omega(|G|)=k$. By \cite[Lemma 1.3]{Keller2} we may assume that $G$ is $\sigma$-reduced, that is, every non-trivial Sylow subgroup of $G$ is isomorphic to a chief factor of $G$. Also write
$n_i=\omega(|G^{(i)}/G^{(i+1)}|)$ for $i\geq 0$, where $G^{(i)}$ is the $i$th term of the derived series of $G$ for $i\geq 0$. As is common, we write $G$, $G'$, $G''$, and $G'''$ for $G^{(i)}$, $i=0,1,2,3$, respectively. Clearly $n_i\leq n$ for all $i$. So we can write $n_{0} + n_{1} + n_{2} = 3n - \lambda$ for some integer $\lambda \geq 0$. We therefore have that $\omega(|G'''|) = k - 3n + \lambda$. Now let $\tau$ be the set of primes less than or equal to $\omega(|G'''|)$, and let $H$ be a Hall $\tau'$-subgroup of $G$ (so that all prime divisors of $|H|$ are greater than $\omega(|G'''|)$). Again by \cite[Lemma 1.3]{Keller2} let $H_0$ be a $\sigma$-reduced subgroup of $H$ such that $\omega(|H_0|)=\omega(|H|)$. Then by 
\cite[Theorem 7(a)]{Keller}  we have that $\omega(|H|)=\omega(|H_0|)\leq 4\omega(H_0)\leq 4\omega(H)\leq 4n$. Thus
\[\rho(n) =k=\omega(|G|)\leq \omega(|H|)+|\tau| \leq 4n + \pi(k-3n + \lambda).\] 
If the right-hand side of this inequality is $\leq k-1$ we have a contradiction. This means we need those values of $\lambda$ such that
\begin{equation}\label{plum}
\pi(k-3n + \lambda) \leq k - (4n +1).
\end{equation}
For a given $n$ and a fixed $k$ it is clear that if  (\ref{plum}) is satisfied at all, then it will only be satisfied for $\lambda \leq L$ for some $L$. Hence, for all/any such $\lambda$ we derive a contradiction, and so we can eliminate the case $\rho(n) = k$.

For larger values of $\lambda$ we use a variation of (\ref{oboe}). By following the proof of the second author \cite{Keller} we have $\rho(n) \leq \pi(\rho(n)) + 4n - \lambda$. Since in these cases we have $\lambda\geq L+1$, we obtain a contradiction if
\begin{equation}\label{mustard}
4n - L \leq k - \pi(k).
\end{equation}
For a given $n$ and a fixed $k$, it is easier to solve (\ref{mustard}) if $L$ is large, that is, if there are many $\lambda$ for which (\ref{plum}) is true.

Let us illustrate this method when $n=9$: we already know that $k\leq 51$. Suppose $k=51$. Then (\ref{plum}) is true for all $0\leq \lambda \leq 22$, and so, for such $\lambda$ we can conclude that $\rho(9) \leq 50$. For larger $\lambda$ we turn to (\ref{mustard}): for $L=22$ the left-hand side is 14, whereas the right-hand side is 36, hence we still can conclude that $\rho(9) \leq 50$.

We continue in this way and can eliminate $\rho(9) \in[46,50]$ but cannot rule out $\rho(9) = 45$. There, (\ref{plum}) is only satisfied for $0\leq \lambda \leq 4$, whence $L=4$, but then (\ref{mustard}) is not satisfied. Hence, the best we can do is to conclude that $\rho(9) \leq 45.$

This approach yields the following improved upper bounds on $\rho(n)$:
\begin{equation}\label{method}
(n, C(n)) = \{(8, 40), (9, 45), (10, 49), (11, 53), (12, 57)\}.
\end{equation}
We find that the worst bounds for $\rho(n)/n$ arise when $n=8$ and $n=9$. The values on (\ref{clarinet}) and (\ref{method}), and a verification on a standard desktop computer of the above method establishes that $\rho(n)\leq 5n$ for all $n\leq 265$. 

When $n\geq 266$, we use the bound on $\pi(n)$ from (\ref{anatra}) in (\ref{oboe}), obtaining $$\rho(n)\leq 1.25506\cdot\frac{\rho(n)}{\log (\rho(n))}+4n.$$
	Now, using the third relation in (\ref{flute}), namely that $\rho(n)\geq 2n$, we have$$\rho(n)\left(1-1.25506\cdot\frac{1}{\log(2n)}\right)\leq\rho(n)\left(1-1.25506\cdot\frac{1}{\log(\rho(n))}\right)\leq 4n.$$For $n=266$ we have$$4\cdot\left(1-1.25506\cdot\frac{1}{\log(2\cdot 266 )}\right)^{-1}=4.9997\ldots<5.$$
	% 4.99973245703
%	In addition, for $n\geq 266$ we have\begin{equation*}\begin{aligned}\rho(n)\left(1-1,25506\cdot\frac{1}{\log(2\cdot 266)}\right)&\leq\rho(n)\left(1-1,25506\cdot\frac{1}{\log(2n)}\right)\\&\leq\rho(n)\left(1-1,25506\cdot\frac{1}{\log(\rho(n))}\right)\\&\leq 4n.\end{aligned}\end{equation*}
As a result, we can conclude that $\rho(n)\leq 5n$ for $n\geq 266.$
%\end{enumerate}
%$$$$We can conclude that $$\rho(n)\leq 5n\quad\forall n\geq 17.$$

We remark in passing that this part of the proof could be substantially improved by using better results for $\pi(n)$. Indeed, Rosser and Schoenfeld \cite{RS} give other bounds that can be re-written into the form of (\ref{bassoon}) as required. However, as noted earlier, the difficult cases come from small values of $n$: here for $n=8,9$.

We note that the improvement for $\rho(9)$, and a quick inductive argument allows one to improve the second bound in (\ref{flute}) to 
$$\rho(n) \leq \frac{n(n+3)}{2} - 9, \quad (n\geq 9).$$

We also remark that we have used the bound $\rho(n)\geq 2n$, which holds for all $n$. We could obtain sharper results, for larger $n$, if we knew
that even $\rho(n)\geq 3n$ for $n\geq 4$. For example, by using this bound, the $4.9997\ldots$ above could be improved to 4.925041. Now note that by
\cite[Remark 1.11]{Keller2} we already know that $\rho(n)\geq 3n$ whenever $n$ is divisible by 4. One gets this by taking direct products
of suitable groups as constructed in \cite[Example 1.9]{Keller2}. To prove that $\rho(n)\geq 3n$ for all $n\geq 4$, all one needs to do is to construct
three groups, using the same method as in \cite[Example 1.9]{Keller2}, to prove that, in the notation of \cite{Keller2}, $\rho_3(5)=15$,
$\rho_3(6)=18$, and $\rho_3(7)=21$, respectively. With these groups, and taking direct products with groups as in \cite[Example 1.9]{Keller2},
all values of $n\geq 4$ can be reached, and so $\rho(n)\geq 3n$ for $n\geq 4$. Since the expected improvements to our bounds in this
paper seem to be minor, we leave the construction of the three above-mentioned groups to the interested reader.\\

%We believe this bound to be true for all $n\geq 3$, and may return to this in future work%.\footnote{\textbf{From TST: Some comment, Thomas, on this? Or leave it as is is?}}
\section{An application to a related problem and a proof of Corollary \ref{bistecca}}\label{foca} 

There is also a related problem for arbitrary finite $G$, that is, not just for the finite solvable groups we have considered hitherto. Whereas for solvable $G$ we have a linear bound of $\rho(n)$ in $n$, for arbitrary $G$ the sharpest result to date is by Hung and Yang \cite{Hung}, who show that $\rho_{g}(n) \leq 210n^{4}$.

We note that, towards the bottom of page 1910 in \cite{Hung}, the proof gives the following bound 
\begin{equation}\label{padella}
\rho_{g}(n) \leq n^{4}\left(28 C(n) + \frac{C(n)^{2}}{n^{2}} + \frac{C(n)}{n^{3}}\right).
\end{equation}
Hung and Yang use a weaker version of our Theorem \ref{gelato} with $\rho(n) \leq 7n$ for all $n$. They also prove that $\rho_{g}(1) \leq 4$, and so, with $n\geq 2$, they can take the term in the parentheses of (\ref{padella}) to be $209.125<210$.

We obtain values for $C(n)$ from (\ref{clarinet}), (\ref{method0}) and (\ref{method}), and use Theorem \ref{gelato}. This shows that the term in the parentheses in (\ref{padella}) is at most $140.41$, which proves Corollary \ref{bistecca}. 

We have not made any other attempt at improving the result. 
A close inspection of the proof in \cite{Hung} shows that, at least when $n$ is large, one should be able to replace the 28 in (\ref{padella}) by 8. This, and the fact that $\lim_{n\rightarrow\infty}C(n) = 4,$ means that one has $\rho_{g}(n) \leq (32+ \epsilon) n^{4}$ for $n$ sufficiently large, and for any positive $\epsilon$.

We conclude by noting that the method used in \S \ref{delfino},  could also be deployed in the paper by the second author \cite{Keller97}. Finally, since specific bounds on $\pi(x)$  are used by Moret\'{o} in the proofs of Theorem 1.1 and 1.2 in \cite{MoretoIMRN}, it is possible that our methods may be of some use there as well.
%: specifically, his Theorem 1.1 and 1.2. Note that bounds on $\pi(x)$ are used in his Lemma 2.2. See also \cite[pp.\ 3381-3382]{MoretoIMRN}.
\section*{Acknowledgements}

We wish to thank Mike Mossinghoff and Alexander Moret\'{o} for helpful discussions on this topic. 
%\textbf{Add others as appropriate.}

 \end{document}